\documentclass[11pt]{amsart}
\usepackage{psfrag,epsf,amsmath,amscd,amsthm,amsfonts,latexsym}
\usepackage[dvips]{graphicx}

\theoremstyle{plain}

\newtheorem{theorem}                 {Theorem}      [section]
\newtheorem{proposition}  [theorem]  {Proposition}

\newtheorem{lemma}        [theorem]  {Lemma}
\theoremstyle{definition}

\newtheorem{example}      [theorem]  {Example}
\newtheorem{remark}       [theorem]  {Remark}

\numberwithin{equation}{section}

\def \r{\mbox{${\mathbb R}$}}
\def \h{\mbox{${\mathbb H}$}}
\def \s{\mbox{${\mathbb S}$}}

\DeclareMathOperator{\Divv}{div_{\tiny\h^2\times\r}}

\DeclareMathOperator{\Div}{div}
\DeclareMathOperator{\R}{R}
\begin{document}

\begin{abstract}
In this article we consider surfaces in the  product space
$\h^2\times \r$ of the hyperbolic plane  $\h^2$ with the real
line. The main results are:  a description of  some geometric properties of minimal graphs;
new examples of complete minimal graphs;
the classification of umbilical surfaces.
\end{abstract}

\title{A note on surfaces in $\h^2\times\r$}

\author{Stefano Montaldo}
\address{Universit\`a degli Studi di Cagliari\\
Dipartimento di Matematica e Informatica\\
Via Ospedale 72\\
09124 Cagliari, Italia}
\email{montaldo@unica.it}

\author{Irene I. Onnis}
\thanks{Work partially supported by GNSAGA and  INdAM, Italy}
\address{Departamento de Matem\'{a}tica, C.P. 6065\\
IMECC, UNICAMP, 13081-970, Campinas, SP\\ Brasil}
\email{onnis@ime.unicamp.br}



\maketitle

\section{Introduction}

In the last decade the study of the geometry of surfaces in the
three-dimensional  Thurston geometries has grown considerably.
One reason is that these spaces can be endowed with a complete metric
with a large isometry group; another, more recent, is the announced 
proof of Thurston geometric conjecture, which ensures the dominant role of this spaces among the three-dimensional geometries.

Leaving aside the space forms $\r^3$, $\s^3$ and $\h^3$, among the left
five Thurston geometries the Heisenberg space is probably the most studied
and the geometry of  surfaces is well understood. In recent years the study of the geometry of surfaces in the two
product spaces $\h^2\times\r$ and $\s^2\times\r$ is growing very rapidly,
and the interest is mainly focused on minimal and constant mean curvature surfaces  \cite{AR,cad,FM,MO,MO2,MR,NR,NR1,R,S,ST}.

The purpose of this paper is first to investigate on some
geometric properties of minimal graphs in $\h^2\times\r$
(Theorem~\ref{pro-geo} and \ref{teo-min})
and to produce some new examples including complete ones.
In the last part (Theorem~\ref{uffa}) we classify the umbilical
surfaces in $\h^2\times\r$ giving their explicit local parametrizations.

For completeness we recall some basic notions on $\h^2\times\r$.
Let $\h^{2}$ be the upper half-plane model $\{(x,y)\in\r^2\,|\,y>0\}$
of the hyperbolic plane endowed with the metric
$g_{\tiny\h}=(dx^2+dy^2)/y^2$, of constant Gauss curvature $-1$.
The space $\h^2$, with the group structure derived by the
composition of proper affine maps, is a Lie group and the
metric $g_{\tiny\h}$ is left invariant. Therefore the
product $\h^2\times\r$ is a Lie group with the left
invariant product metric
$$
g=\frac{dx^2+dy^2}{y^2}+dz^2.
$$
With respect to the metric $g$ an orthonormal basis of
left invariant vector fields  is
\begin{equation}\label{campi}
E_1=y\frac{\partial}{\partial x},\qquad E_2=y\frac{\partial}{\partial y},\qquad
E_3=\frac{\partial}{\partial z},
\end{equation}
and the non zero components of the Christoffel symbols are:
\begin{equation}\label{chri}
\Gamma_{12}^1= \Gamma_{21}^1=\Gamma_{22}^2=-\frac{1}{y},\qquad \Gamma_{11}^2=\frac{1}{y}.
\end{equation}

\section{Minimal graphs}
The natural parametrization of a graph
$\mathcal{M}$ in $\h^2\times\r$ is
$$
\phi(x,y)=(x,y,f(x,y)),\qquad (x,y) \in \Omega,
$$
where the domain $\Omega\subseteq \h^2$ is relatively compact, with a
differentiable boundary and  $f:\Omega\to\r$ is a $C^2$-function.
The unit normal $\xi$ to $\mathcal{M}$ is given by
\begin{equation}\label{xi}
\xi(x,y)=-\frac{f_x}{w\,y}\,E_1-\frac{f_y}{w\,y}\,E_2+\frac{1}{w\,y^2}\,E_3,
\end{equation}
where $w=\tfrac{1}{y^2}\sqrt{y^2(f_x^2+f_y^2)+1}$.
The coefficients of the induced metric $h=\phi^{\ast}g$ are
$$
E=f_x^2+\frac{1}{y^2},\qquad F=f_x f_y,\qquad
G=f_y^2+\frac{1}{y^2},
$$
while the coefficients of the second fundamental form are given by
\begin{equation}\label{lmn}
L=\frac{y f_{xx}-f_y}{w\, y^3},\qquad M=\frac{y f_{xy}+f_x}{w
\,y^3},\qquad N=\frac{y f_{yy}+f_y}{w\, y^3}.
\end{equation}
The mean curvature function is then 
\begin{equation}\label{acca}
H=\frac{y^2}{2}\,\Div\Big(\frac{\nabla f}{\sqrt{1+y^2\,|\nabla
f|^2}}\Big)=\frac{1}{2}\,\Div_{\tiny\h^2}\Big(\frac{\nabla f
}{w}\Big),
\end{equation}
where $\nabla$ and $\Div$ stand for the Euclidean gradient and the Euclidean 
divergence, while $\Div_{\tiny\h^2}$ is the divergence in $(\h^2,g_{\tiny\h})$. 
The equation $H=0$ is called the {\it minimal
surfaces equation} in $\h^2\times\r$, and can be also written as
\begin{equation}\label{eq1}
(1+y^2f_y^2)\,f_{xx}-y\,(f_x^2+f_y^2)\,f_y-2y^2\,f_x\,f_y\,f_{xy}+(1+y^2f_x^2)\,f_{yy}=0.
\end{equation}
This equation was first founded by  B.~Nelli and H.~Rosenberg, in \cite{NR}, 
where they showed that in $\h^2\times\r$
there exist minimal surfaces of Catenoid-type, Helicoid-type and Scherk-type.
Moreover, they proved that the Bernstein's theorem fails, that is there exist
complete minimal graphs in $\h^2\times\r$ of rank different from zero.

The first geometric property of minimal graph is that,
as in the Euclidean case (see \cite{GDG}), solutions of \eqref{eq1} define graphs of ``minimal''
area. 

\begin{theorem}\label{teo-min}
If $f$ satisfies the minimal surfaces equation \eqref{eq1} in $\Omega$ and $f$ extends
continuously to  $\overline{\Omega}$, then the area of the surface $\mathcal{M}$, defined by $f$,
is less than or equal to  the area of any other surface $\mathcal{\widetilde{M}}$ defined by a
function $\widetilde{f}$ in $\Omega$ having the same values of $f$ on  $\partial\Omega$. Moreover,
equality holds if and only if $f$ and $\widetilde{f}$ coincide on $\Omega$.
\end{theorem}
\begin{proof}
In the domain $\Omega\times\r$ of $\h^2\times\r$, consider the
unit vector field $V(x,y,z)$ given by
$$
V(x,y,z):=-\frac{f_x}{w\,y}\,E_1-\frac{f_y}{w\,y}\,E_2+\frac{1}{w\,y^2}\,E_3.
$$
Writing
$V=V^i\,(\partial/\partial x_i)$ and denoting by $\Divv$
 the divergence of $\h^2\times\r$, we have
$$
\Divv V=-y^2\,\Div\Big(\frac{\nabla f}{\sqrt{1+y^2\,|\nabla
f|^2}}\Big).
$$
Since $f(x,y)$ satisfies Equation~{\eqref{eq1}},
it follows that
$$
\Divv V\equiv 0 \qquad\text{on}\quad
\Omega\times\r.
$$
The surfaces $\mathcal{M}$ and
$\mathcal{\widetilde{M}}$ have the same boundary, and therefore
$\mathcal{M}\cup\mathcal{\widetilde{M}}$ is an oriented boundary of
an open set $\Theta$ in $\Omega\times\r$. Denoting by
$\eta$ the unit normal corresponding to the positive orientation
on $\mathcal{M}\cup\mathcal{\widetilde{M}}$ and using the
Divergence Theorem, we have:
\begin{equation}\label{dive}
  0=\int_\Theta \Divv
  V=\int_{\mathcal{M}\cup\mathcal{\widetilde{M}}}g(V,\eta)\,dA\,.
\end{equation}
From the definition of the vector $V$ and from {\eqref{xi}}, it
follows that
$$
V\equiv \eta \qquad\text{on}\quad \mathcal{M},
$$
hence, from {\eqref{dive}}, and since $V$ and $\eta$ are both unit vector fields, it results that
\begin{align*}
A(\mathcal{M})=\int_{\mathcal{M}}g(V,\eta)\,dA=\int_{\mathcal{\widetilde{M}}}g(V,-\eta)\,dA
\leq \int_{\mathcal{\widetilde{M}}}dA=A(\mathcal{\widetilde{M}}).
\end{align*}
Furthermore, equality holds if and only if  $g(V,-\eta)=1$, that is, if and only if
on $\Omega$ $\widetilde{f}_x=f_x$ and $\widetilde{f}_y=f_y$.
Finally, since
$f_{|{\partial \Omega}}=\widetilde{f}_{|{\partial \Omega}}$, we must have
that $\widetilde{f}(x,y)=f(x,y)$ for all $(x,y)\in\Omega$.
\end{proof}

In the following we show some solutions
of the Equation {\eqref{eq1}}.
\begin{example}\label{planoplano}
If a solution of {\eqref{eq1}} has the form $f(x,y)=\varphi(x)$,
we have that ${\varphi}''(x)=0$ and, then, $f(x,y)=ax+b$, with $a,b\in\r$.
These are the only minimal planes in $\h^2\times\r$ that can be described as graphs.

If now we look for solutions of type $f(x,y)=\psi (y)$,  then
Equation~{\eqref{eq1}} assumes the form
${\psi}''(y)-y\, {\psi}'(y)^3=0$, and, by integration, we get
$$
\psi (y)=\arcsin (a\, y)+b,\qquad
0<y\leq 1/a\,,\quad a,b\in\r,\quad a>0.
$$
\end{example}

\begin{example}\label{funnel}
As for the Euclidean space, we can find interesting examples of minimal graphs seeking for radial
solutions of {\eqref{eq1}}  of type $f(x,y)=h(x^2+y^2)$. In this case, it results that
$z\,h''(z)+h'(z)=0$, thus the desired function is
$$
f(x,y)=a\,\ln(x^2+y^2)+b,\qquad a,b\in\r.
$$
This surface, called the {\it funnel surface}, defines a complete minimal
graph.
We observe that the Gauss map of this surface is of rank $1$.
On the right hand side of Figure~\ref{gauss}
there is a plot of the image,  under the Gauss
map, of the funnel surface, which is plotted on the left hand side.
\begin{figure}[h]
\begin{center}
\begin{minipage}[c]{1.1in}
\psfragscanon \psfrag{x}[r][t][1]{$x$}
        \psfrag{y}[t][1]{$y$}
         \psfrag{z}[l][1]{$z$}
\includegraphics[width=1.1\textwidth]{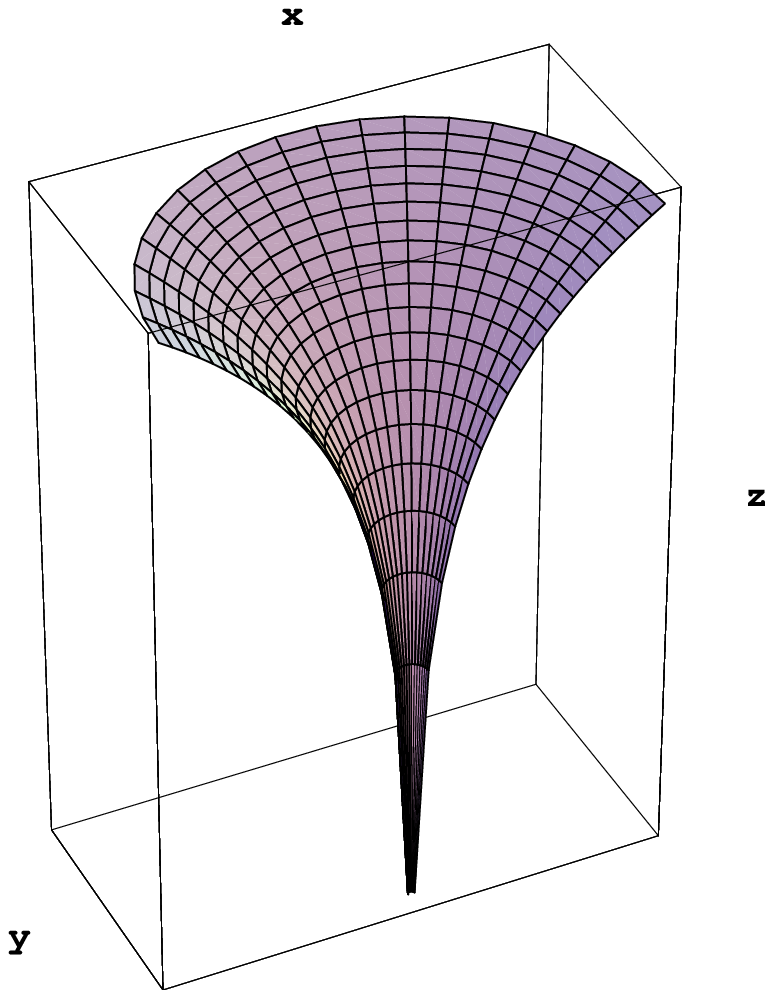}
\end{minipage}
\hspace{20mm}
\begin{minipage}[c]{1.1in}
   \includegraphics[width=1.3\textwidth]{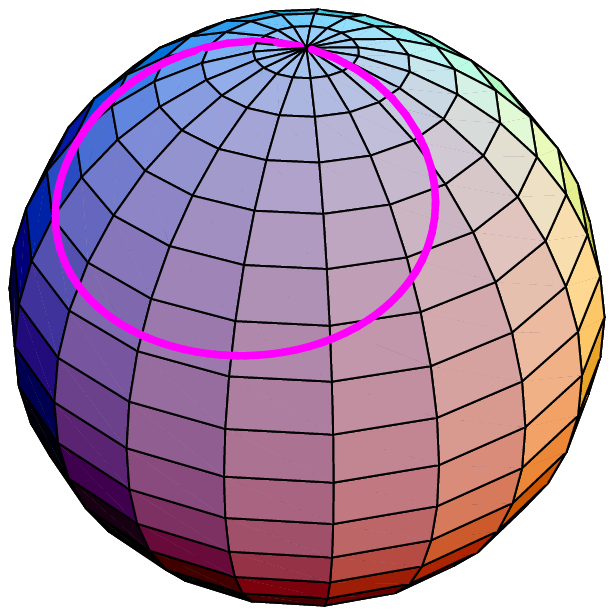}
\end{minipage}
\end{center}
\caption{Complete minimal graph of rank $1$ (left) and the image of its Gauss map
(right).}\label{gauss}
\end{figure}
\end{example}

\begin{example}\label{exnuova}
Let $f(x,y)$ be a solution of the minimal surfaces equation of type
$$
f(x,y)=\frac{a(x)}{x^2+y^2},
$$
where $a(x)$ is a real function. Then, {\eqref{eq1}} gives
$$
[(x^2+y^2)^4+4\,y^4\,a(x)^2]\,a''(x)-4(x^2+y^2)^3\,[x\,a'(x)-a(x)]=0,
$$
of which a solution is $a(x)=c\, x$, with $c\in\r$. The corresponding
minimal function is
\begin{equation}\label{eq-min-r2}
f(x,y)=\frac{c\,x}{x^2+y^2},\qquad c\in\r,
\end{equation}
which produces the minimal graph plotted in Figure~\ref{nuovaa} (left). 
We observe that the Gauss map of this complete graph is of rank $2$.

This example can be generalized considering, for a given real function $h$, a 
solution of  \eqref{eq1}  of type
$f(x,y)=h\big(\tfrac{x}{x^2+y^2}\big)$ or of type
$f(x,y)=h\big(\tfrac{y}{x^2+y^2}\big)$.
In the first case we essentially find (up to translations) the example given
by \eqref{eq-min-r2}. In the second case it  results that
$h''(z)-z\,h'(z)^3=0$ and, therefore,
$$
h(z)=\arcsin (a z)+b,\qquad a,b\in\r.
$$
The corresponding minimal function does not define a complete graph. 
A plot of this surface is given in Figure~\ref{nuovaa} (right).

\begin{figure}[h]
\begin{center}
\begin{minipage}[b]{1.6in}
\centering \psfragscanon \psfrag{x}[l][1]{$x$}
        \psfrag{y}[b][r][1]{$y$}
         \psfrag{z}[b][1]{$z$}
\includegraphics[width=1.\textwidth]{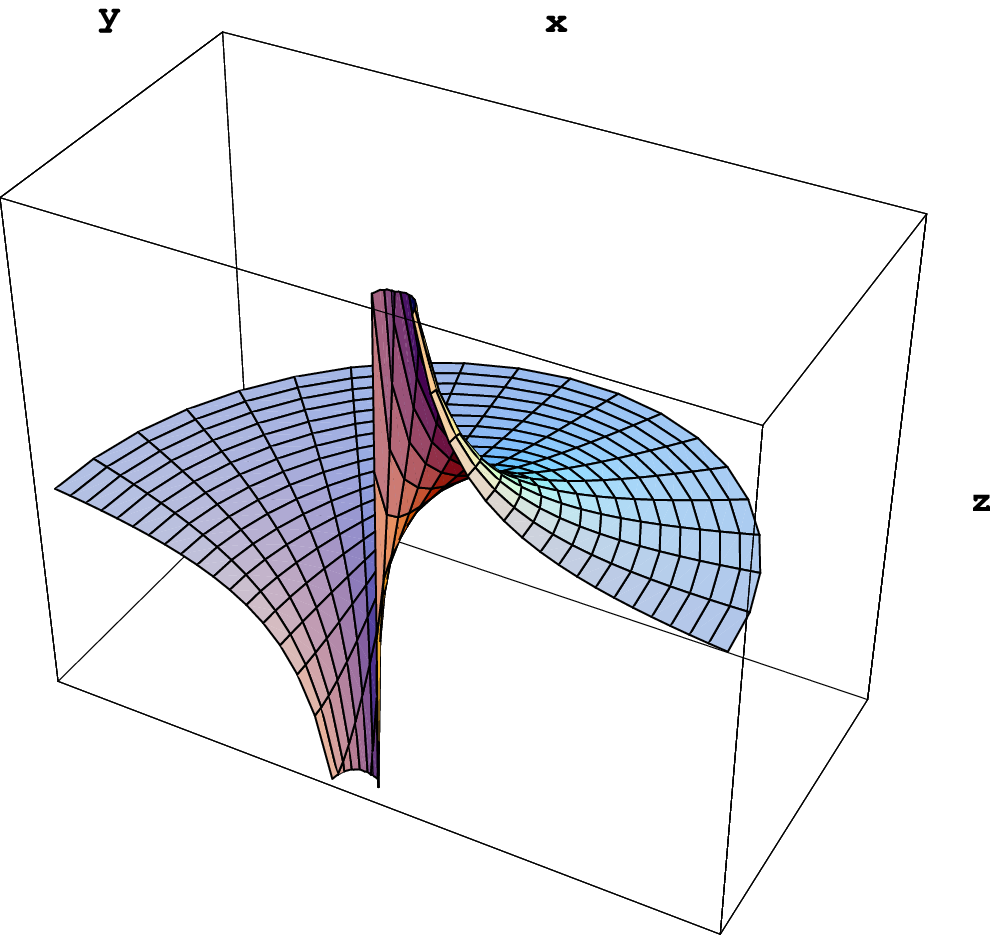}
\end{minipage}
\hspace{15mm}
\begin{minipage}[b]{1.6in}
    \centering \psfragscanon \psfrag{x}[l][1]{$x$}
        \psfrag{y}[b][r][1]{$y$}
         \psfrag{z}[b][1]{$z$}
\includegraphics[width=1.\textwidth]{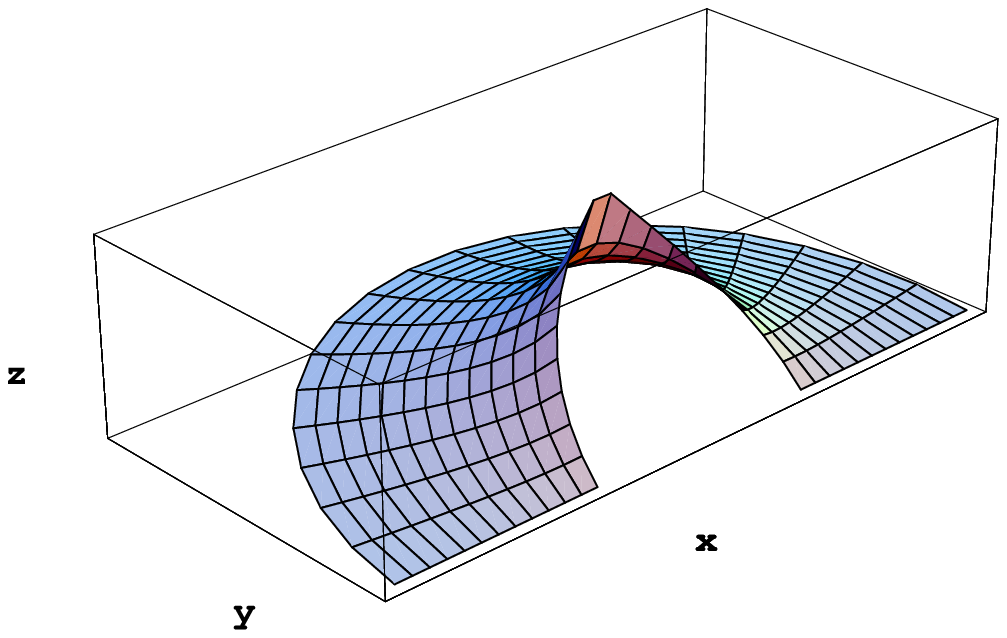}
\end{minipage}
\end{center}
\protect\caption{Minimal graphs of rank $2$: a complete one (left)  
and a non complete one (right).}\label{nuovaa}
\end{figure}
\end{example}

\section{Minimality and harmonicity}
In this section we study the relations between the minimality of a
surface $\mathcal{M}$ in $\h^2\times\r$, defined as the graph of a differentiable
function $f$, and the harmonicity of $f$. A first property about
minimal graphs in $\h^2\times\r$ is given by the following:

\begin{proposition}
Let $\mathcal{M}\subset \h^2\times\r$ be a graph of a $C^2$-function $f$, 
defined in a domain $\Omega$ of $\h^2$.
If $\mathcal{M}$ is a minimal graph then
$f$ is harmonic with respect to the induced metric $h=\phi^{\ast}g$, where
$\phi(x,y)=(x,y,f(x,y)), (x,y)\in\Omega$, is a global parametrization of
$\mathcal{M}$.
\end{proposition}
\begin{proof}
We use (see, for example, \cite{ES}) the fact that
$\mathcal{M}$ is a minimal surface if and only if
 $\phi:(\Omega,h)\to\h^2\times\r$ is harmonic, that is $\phi$ satisfies the system
$$
\Delta_{h}\phi^{\alpha}+h^{ij}\,\Gamma_{\beta\gamma}^{\alpha}\,\frac{\partial
\phi^{\beta}}{\partial x^i} \frac{\partial \phi^{\gamma}}{\partial
x^j}=0, \quad \alpha=1,2,3,
$$
where $\Delta_{h}$ is the Beltrami-Laplace operator with respect to $h$. Since
$\Gamma_{\beta\gamma}^{3}=0$, for all $\beta,\gamma\in\{1,2,3\}$, (see \eqref{chri}), it follows
that $\Delta_h f=0$.
\end{proof}

As an immediate consequence of the above proposition we can prove that
 there exist no compact minimal surfaces  without boundary in the product $\h^2\times\r$.

The next result gives a link between the harmonicity of $f$
and the geometry of the level curves.

\begin{proposition}\label{pro-geo}
Let $\mathcal{M}\subset \h^2\times\r$ be a minimal surface defined as the graph of a non constant
differentiable function $f:\Omega\subseteq\h^2\to\r$. Then, the level curves  of
$f$ are pre-geodesics of $\h^2$ if and only if $f$ is harmonic with respect to the flat
Laplacian.
\end{proposition}
\begin{proof}
Let $f$ be a solution of {\eqref{eq1}} and let  $\gamma(t)=(x(t),y(t))$ be the 
parametrization of a level curve of $f$. If the function $x(t)$ is constant, it 
results that $f_y=0$ and thus, using  Example
\ref{planoplano}, the surface $\mathcal{M}$ is a piece of the plane $z=ax+b$,
$a,b\in\r$.

Thus we can assume that there exists a point $t_0$ such that $x'(t)\neq 0$ in a 
neighborhood of $t_0$. Therefore, we can
parametrize $\gamma$ as $\gamma(x)=(x,y(x))$, with  $y(x)>0$,
in a neighborhood of  $\gamma(t_0)$. It follows that
$$
\gamma'(x)=\frac{E_1}{y(x)}+\frac{y'(x)}{y(x)}\, E_2,
$$
and
$$
\nabla_{\gamma'}\gamma'=-\frac{2y'(x)}{y(x)^2}\,E_1+
\Big(\frac{y(x)y''(x)-y'(x)^2+1}{y(x)^2}\Big)\,E_2.
$$
The geodesic curvature $k_g$ (in $\h^2$) of $\gamma$ is then
$$
k_g=\frac{g(\nabla_{\gamma'}\gamma',J\gamma')}{\|\gamma'\|^3}
=\frac{y(x)y''(x)+y'(x)^2+1}{(1+y'(x)^2)^{3/2}}.
$$
Using 
$$
 y'(x)=-\frac{f_x(x,y(x))}{f_y(x,y(x))}
$$
and
$$
y''(x)=-\frac{f_{xx}+2f_{xy}y'+f_{yy}(y')^2}{f_y},
$$
the minimal Equation~\eqref{eq1} can be written as
$$
\Delta f=y
(x)\,k_g\,|\nabla f|_{\tiny \r^2}^3,
$$
which completes the proof.
\end{proof}

\section{Umbilical surfaces of $\h^2\times\r$}
We start this section studying the totally geodesic surfaces of
$\h^2\times\r$. For this we need the following lemma.
\begin{lemma}\label{lemver}
Let $\mathcal{M}$ be a regular surface in $\h^2\times\r$. Then, there
exists an open dense set  in $\mathcal{M}$, of which the connected
components admit one of the following parametrizations:
$$\begin{aligned}
X(u,v)&=(u,v,f(u,v)),\;\quad v>0,\\
Y(u,v)&=(u,a(u),v),\;\qquad a(u)>0,\\
Z(u,v)&=(c,u,v),\;\qquad u>0.
\end{aligned}
$$
\end{lemma}
\begin{proof}
A detailed proof can be found in \cite{O}.
\end{proof}
The  parametrizations $Y(u,v)$ and $Z(u,v)$ of Lemma~\ref{lemver} define surfaces that we shall call
{\it vertical surfaces}.

\begin{theorem}\label{pintus} The totally geodesic surfaces of $\h^2\times\r$
are the horizontal planes $z=c$, $c\in\r$, and the vertical
cylinders over the geodesics  of $\h^2$.
\end{theorem}
\begin{proof} We start proving that the only totally geodesic graphs  are the horizontal planes
$z=c$, $c\in\r$. Let $\mathcal{M}$ be a totally geodesic surface
defined as the graph of a differentiable function $f$. From
{\eqref{lmn}}, it follows that
\begin{equation}\label{sal}
\left\{
\begin{aligned}
f_{xx}&=\frac{f_y}{y},\\
f_{yy}&=-\frac{f_y}{y},\\
f_{xy}&=-\frac{f_x}{y}.
\end{aligned}\right.
\end{equation}
First, observe that $f(x,y)=c$, $c\in\r$, satisfies System~\eqref{sal}  and, 
therefore, defines a totally geodesic surface.
Then, since $f_x=0$ if and only if $f_y=0$, there exist no totally geodesic
 graphs defined by a (non constant) function $f$
that depends only of one variable.
Thus, we can suppose that
$f_x\neq 0$ and $f_y\neq 0$. From the second and third equations of
{\eqref{sal}}, we find that there exist two functions $p(x)$ and $q(x)$
such that $f_y={p(x)}/{y}$ and  $f_x={q(x)}/{y}$. Now, replacing in the first 
equation of {\eqref{sal}}, we have the contradiction
$$
y\,q'(x)-p(x)=0.
$$

To complete the proof, observe that the vertical
cylinders over the geodesics of $\h^2$ are totally geodesic surfaces of
$\h^2\times\r$. In the following, we prove that these cylinders
are the only vertical surfaces that are totally geodesic in
$\h^2\times\r$. Let $\mathcal{M}$ be a vertical surface. From
Lemma~\ref{lemver}, it follows that either $\mathcal{M}$ is the
totally geodesic plane $x=c$, with $y>0$, or it is parametrized by
\begin{equation}\label{1fo}
X(u,v)=(u,a(u),v),\qquad a(u)>0.
\end{equation}
The unit normal to the surface $\mathcal{M}$, defined by \eqref{1fo},
is
\begin{equation}\label{xi1}
\xi=\frac{{a}'}{\sqrt{1+({a'})^2}}\,E_1-\frac{1}{\sqrt{1+({a'})^2}}\,E_2.
\end{equation}
It is then straightforward to compute that
$$
\nabla_{X_u}\xi=\bigg[\Big(\frac{{a'}}{\sqrt{1+({a'})^2}}\Big)_u+
\frac{1}{a\,\sqrt{1+({a'})^2}}\bigg]\,E_1+
\bigg[\frac{{a'}}{a\,\sqrt{1+({a'})^2}}-
\Big(\frac{1}{\sqrt{1+({a'})^2}}\Big)_u\bigg]\,E_2,
$$
and
$$
\nabla_{X_v}\xi=0,
$$
hence $M=N=0$. Consequently, $\mathcal{M}$ is totally geodesic if and
only if
$$
L=-g(\nabla_{X_u}\xi,X_u)=0,
$$
that is, the function $a(u)$ satisfies the following ODE:
$$
a\,{a''}+({a'})^2+1=0.
$$
This implies that
$$
a(u)=\sqrt{-u^2+2\,c_1\,u+c_2},\qquad
c_1,c_2\in\r\quad\text{with}\quad c_1^2+c_2>0,
$$
and the curve $(u,a(u))$ is the geodesic of $\h^2$ given by the upper
semi-circle with center at $(c_1,0)$ and radius $\sqrt{c_1^2+c_2}$.
This completes the proof.
\end{proof}

\begin{remark}
From the proof of Theorem~\ref{pintus} it  follows that, if
$\mathcal{M}$ is a vertical surface, then $F=M=N=0$ and so the
mean curvature is given by
$H={L}/{2 E}.$ Thus a vertical surface is  minimal
if and only if it is totally geodesic.
\end{remark}

We are now ready to state the main result of this section.

\begin{theorem}\label{uffa}
  The umbilical surfaces of $\h^2\times\r$ are:
\begin{enumerate}
\item[i)]  the totally geodesic surfaces given in Theorem~\ref{pintus};
\item[ii)] the surface given as the graph of the function
$$
f(x,y)=\arctan\Big(\frac{\lambda(x,y)}{\sqrt{j-\lambda(x,y)^2}}\Big)+c,
\qquad c\in\r$$
where
$$
\lambda(x,y)=\frac{1}{y}\,\Big[\frac{c_1}{2}
(x^2+y^2)+c_2\,x-c_3\Big],\qquad c_1,c_2,c_3\in\r,
$$
and
$$
j=1-c_2^2-2\,c_1\,c_3>0.
$$
\end{enumerate}

\end{theorem}
\begin{proof}
Let $\mathcal{M}$ be an umbilical surface of $\h^2\times\r$
parametrized  by
$$
\phi(x,y)=(x,y,f(x,y)), \quad (x,y)\in\Omega\subseteq\h^2.
$$
Let
$p\in \mathcal{M}$ and let $\xi$ be an unit normal vector field defined in
some neighborhood $U$ of $p$. Since $\mathcal{M}$ is umbilical, by definition,
there exists a function $\lambda: U\rightarrow \r$ such that the shape operator $A$ satisfies
$A_{\xi}=\lambda I$ in $U$.
The expression of $A_{\xi}$ with respect to the coordinates basis $\{\phi_x,\phi_y\}$ is
$$
A_\xi (p)=\begin{pmatrix}
\Big(\displaystyle{\frac{f_x}{w}\Big)_x-\frac{f_y}{y\,w}} & \qquad y\Big(\displaystyle{\frac{f_x}{y\, w}\Big)_y }\\
&\\
\Big(\displaystyle{\frac{f_y}{w}\Big)_x+\frac{f_x}{y\,w}} &\qquad
y\Big(\displaystyle{\frac{f_y}{y\,
  w}}\Big)_y
\end{pmatrix}.
$$
Thus  $\mathcal{M}$ is umbilical if and only if
\begin{equation}\label{matriz}
\left\{
\begin{aligned}
\Big(\frac{f_x}{w\,y}\Big)_y &=0,\\
\Big(\frac{f_y}{w}\Big)_x &=-\frac{f_x}{ w\,y},\\
\Big(\frac{f_x}{w}\Big)_x &=\Big(\frac{f_y}{w}\Big)_y.
\end{aligned}\right.
\end{equation}
The first and second equations of \eqref{matriz} imply that there exist 
two functions $a(x)$ and $b(y)$ such that
$$
\frac{f_x}{y w}=a(x),\qquad \frac{f_y}{w}=-\int a(x) dw + b(y).
$$
Then, from the third equation of \eqref{matriz}, we conclude that
$$
a(x)=c_1 x + c_2,\qquad b(y)=\frac{c_1}{2} y^2 +c_3,
$$
where $c_1,c_2,c_3\in\r$. Thus ${f_y}/{w}=\frac{c_1}{2}(y^2-x^2)-c_2 x +c_3$, which implies that
\begin{equation}\label{eq-lambda}
\lambda(x,y)=y\Big(\frac{f_y}{y w}\Big)_y= \frac{1}{y}\Big[ \frac{c_1}{2}(x^2+y^2)+c_2 x -
c_3\Big].
\end{equation}
From the Codazzi's equation for umbilical surfaces (see, for example, \cite{Daj})
$$
(\R(\phi_x,\phi_y)\xi)^\top=
\lambda_y \phi_x- \lambda_x \phi_y,
$$
and using
$$
(\R(\phi_x,\phi_y)\xi)^\top=\frac{f_y}{w y^3} E_1- \frac{f_x}{w y^3} E_2,
$$
it follows that
\begin{equation}\label{ung}
\left\{
\begin{aligned}
\lambda_x&=\frac{f_x}{w\,y^2},\\
\lambda_y&=\frac{f_y}{w \,y^2}.
\end{aligned}
\right.
\end{equation}
Now, using the identities 
\begin{align*}\Big(\frac{1}{w\,y^2}\Big)_y&=-y\Big[f_x\,\Big(\frac{f_x}{w\,y}\Big)_y+f_y\,\Big(\frac{f_y}{w\,y}\Big)_y\Big]
\\\Big(\frac{1}{w\,y^2}\Big)_x&=-\Big[f_x\,\Big(\frac{f_x}{w}\Big)_x+f_y\,\Big(\frac{f_y}{w}\Big)_x\Big],
\end{align*}
and \eqref{eq-lambda}, a simple calculation gives
\begin{equation}\label{eq-simple}
y^2\,w=\frac{1}{\sqrt{j-\lambda^2}},
\end{equation}
where $j(x)=1-c_2^2-2\,c_1\,c_3>0$.
Substituting \eqref{eq-simple} in the first
equation of {\eqref{ung}} we have
$$
f(x,y)=\int\frac{\lambda_x}{\sqrt{j-\lambda^2}}\,dx=
\arctan\Big(\frac{\lambda}{\sqrt{j-\lambda^2}}\Big)+h(y),
$$
for a certain function $h(y)$.
From the second equation of {\eqref{ung}}, we conclude that $h(y)$ is constant.

To complete the proof, we must study the case when $\mathcal{M}$ is
an umbilical vertical surface. From  Lemma~\ref{lemver}, it follows
that either $\mathcal{M}$ is the totally geodesic plane $x=c$, or it is
given by:
$$
X(u,v)=(u,a(u),v),\qquad a(u)>0.
$$
In the last case, we have $A_{\xi}X_v=-\nabla_{X_v}\xi=0$ and, therefore, $\lambda=0$. 
This conclude the proof of the theorem.
\end{proof}

\begin{figure}[h]
\centering
\psfragscanon
\psfrag{x}[c][c][1][0]{$x$}
\psfrag{y}[1][0]{$y$}
\psfrag{z}[b][l][1][0]{$z$}
\includegraphics[width=0.28\textwidth]{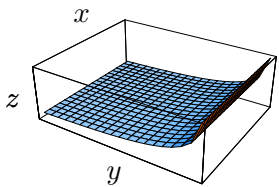}
\caption{\small Umbilical surface of
$\h^2\times\r$.}\label{ombelicale}
\end{figure}

As an example of umbilical surfaces, if $c_1=c_2=0$, we have the graph given by
$$
f(x,y)= -\arcsin\big(\frac{c_3}{y}\big)+c,\qquad y\geq
|c_3|>0.
$$
In  Figure~\ref{ombelicale} there is a plot of this
surface for $c_3=-1$.

\noindent{\bf Acknowledgements}. The authors wish to thank Francesco Mercuri for valuable
conversations during the preparation of this paper.

\end{document}